\documentclass[pdflatex,sn-mathphys-num]{sn-jnl}

\usepackage{graphicx}%
\usepackage{multirow}%
\usepackage{amsmath,amssymb,amsfonts}%
\usepackage{amsthm}%
\usepackage{mathrsfs}%
\usepackage[title]{appendix}%
\usepackage{xcolor}%
\usepackage{textcomp}%
\usepackage{manyfoot}%
\usepackage{booktabs}%
\usepackage{algorithm}%
\usepackage{algorithmicx}%
\usepackage{algpseudocode}%
\usepackage{listings}%

\theoremstyle{thmstyleone}
\newtheorem{theorem}{Theorem}

\theoremstyle{thmstyletwo}
\newtheorem{example}{Example}
\newtheorem{remark}{Remark}

\theoremstyle{thmstylethree}
\newtheorem{definition}{Definition}

\begin{document}

\title[Revisiting Korovkin-type Theorems]{Revisiting Korovkin-type Theorems in Banach Function Spaces}


\author[1]{\fnm{V.B.} \sur{Kiran Kumar}}\email{vbk@cusat.ac.in}

\author*[2]{\fnm{P.C.} \sur{Vinaya}}\email{vinayapc01gmail.com}
\equalcont{These authors contributed equally to this work.}

\affil*[1]{\orgdiv{Department of Mathematics}, \orgname{Cochin University of Science and Technology}, \orgaddress{\street{Kalamassery}, \city{Cochin}, \postcode{682022}, \state{Kerala}, \country{India}}}

\affil[2]{\orgdiv{Department of Mathematics}, \orgname{Cochin University of Science and Technology}, \orgaddress{\street{Kalamassery}, \city{Cochin}, \postcode{682022}, \state{Kerala}, \country{India}}}


\abstract{	This article delves into Korovkin-type theorems in Banach function spaces, as established by Yusuf Zeren et al. (2022). We prove that in this theorem, the positivity of the operators is not a necessary requirement and provide example of a non positive operator where it is applicable. Under the assumption of positivity, we establish an operator version of the result. Additionally, we derive a quantitative form of the result using the modulus of continuity. We apply the result to examples such as Lebesgue space, Weighted Lebesgue space, Grand Lebesgue space, etc. Furthermore, we present numerical illustrations for specific cases.}




\maketitle

\section{Introduction and Preliminaries}\label{sec1}

In $1953,$ P.P. Korovkin introduced an approximation theorem that gives a simple and effective criterion to decide whether a sequence of positive linear operators on $C[0,1]$ constitutes an approximation process \cite{Kor}. This theorem has played a unifying role in various approximation processes and has found applications across different fields in mathematical analysis. Since then, several variations of this theorem have been formulated by various mathematicians, enriching the field with diverse perspectives and applications. A list of relevant results can be found in \cite{korovkin} and recent results can be seen in \cite{altomare1}, \cite{altomare2}. We state the first and second Korovkin theorem.\par
We denote $C[0,1]$ as the space of all continuous functions defined on $[0,1]$ endowed with the sup-norm, $\|f\|_{\infty}=sup\{|f(x)|:x\in [0,1]\}.$
\begin{theorem}{\cite{korovkin}} \label{kor1}
	Let $ \{L_n\}_{n\in\mathbb{N}}$ be a sequence of positive linear operators from $C[0,1]$ to $F[0,1]$ satisfying
	$\lim\limits_{n\rightarrow\infty}L_n(g)=g\ in\ C[0,1]\ for\ all\ g\in\{1,t,t^2\}.$
	Then, $\lim\limits_{n\rightarrow\infty}L_n(f)=f\ in\ C[0,1]\ for\ all\ f\in C[0,1].$
\end{theorem}
Let $C_{2\pi}(\mathbb{R})$ be the space of all continuous and $2\pi-$periodic functions on $\mathbb{R}$ endowed with sup-norm.
\begin{theorem}{\cite{korovkin}}\label{kor2}
	Let $ \{L_n\}_{n\in\mathbb{N}}$  be a sequence of positive linear operators from $C_{2\pi}(\mathbb{R})$ to $F(\mathbb{R})$ satisfying
	$\lim\limits_{n\rightarrow\infty}L_n(g)=g\ in\ C_{2\pi}(\mathbb{R})\ for\ all\ g\in\{1, \sin x, \cos x\}.$
	Then, $
	\lim\limits_{n\rightarrow\infty}L_n(f)=f\ in\ C_{2\pi}(\mathbb{R})\ for\ all\ f\in C_{2\pi}(\mathbb{R}).$
\end{theorem}
In 1968, O. Shisha and B. Mond actively demonstrated quantitative versions of Theorems \ref{kor1} and \ref{kor2}, as outlined in Section $2.3$. In this context, the convergence of a sequence of positive linear operators to the identity operator is achieved by establishing estimates in terms of the convergence on test functions and the modulus of continuity of the function \cite{shisha}.\par 
In $2022$, Yusuf Zeren et al. derived a Korovkin-type theorem within the framework of Banach function spaces \cite{recent}. Their approach involved considering the subspace of the Banach function space defined using the shift operator. By demonstrating the density of , the set of infinitely differentiable functions within this subspace $C_0^{\infty}[0,1]$, they successfully established the Korovkin-type theorem.

They showcased the applications of this result in various concrete examples, including Lebesgue space, Grand Lebesgue space, Morrey space, and their weighted spaces, as well as Weak Lebesgue space, for the Kantorovich polynomials. We give some preliminary definitions and results here . \par 
Let $(A,\mathcal{S},\mu)$ denote a measurable space, where $\mathcal{S}$ is the algebra of the measurable sub sets of $A.$ Let $\mathcal{M}$ be the set of measurable functions on $A$, \textit{$\mathcal{M}^+$} denotes a set of non-negative measurable functions on $A$.
\begin{definition}\cite{recent}
	A mapping $\rho:$ \textit{$\mathcal{M}^+$} $\rightarrow [0,+\infty] $ is called a function norm if the following axioms hold for all $f,g,f_n \in $\textit{$\mathcal{M}^+$}, $a\geq 0$, $E \in \mathcal{S}$:
	\begin{enumerate}
		\item{$\rho(f)=0 \iff f=0, \mu\ a.e.,$
			$\rho(af)=a\rho(f),$ and
			$\rho(f+g)\leq \rho(f)+\rho(g)$};
		\item{$g\leq f,\ \mu\ a.e.\implies\ \rho(g)\leq \rho(f)$ };
		\item{$f_n\uparrow f$, $\mu\ a.e.\ \implies \rho(f_n)\uparrow \rho(f)$};
		\item{$\mu(E)<+\infty\ \implies\ \rho(\chi_E)<+\infty$};
		\item{$\mu(E)<+\infty\ \implies $ there exists $ C_E>0$ such that $\int_{E} fd\mu\leq C_E\rho(f)$, where $C_E$ depends on $E$, $\rho$ and does not depend on $f$.}
	\end{enumerate}
\end{definition}

A Banach function space $X$ generated by $\rho$  is a Banach space of functions $f\in \mathcal{M}$, equipped with the norm 
$		\|f\|_X=\rho(|f|).$ In the following sections, $X$ denotes the Banach function space generated by a function norm $\rho$.
\begin{definition}\cite{recent}
	The function $f\in X$ has an absolutely continuous norm if
	\begin{equation}
		\|f\chi_{E_n}\|_X\rightarrow 0,\ n\rightarrow\infty
	\end{equation}
	for every sequence $E_n\in \Sigma$ satisfying the condition $E_n\rightarrow \emptyset\ \mu\ a.e.$
\end{definition}
Denote 
$X_a=\{f\in X: f$ has an absolutely continuous norm$\}$.	Let $A=[0,1]$, $\mu$ be the Lebesgue measure on the real line, $\mathcal{S}$ be the Borel $\sigma-$algebra. Let $X$ be a Banach function space  with norm $\|.\|_X,$ and $\delta>0$. Consider the shift operator  on $X$, $T_{\delta}$ defined by
\[ T_{\delta}f(x)=
\begin{cases}
	f(x+\delta)\ &if\ x+\delta\in [0,1]\\
	0\ &if\ x+\delta\not\in [0,1]
\end{cases}
\]
We denote $X^S$, as the closure in $X$ of a linear manifold of functions $\{f\in X :
\lim\limits_{\delta\rightarrow 0}\|T_{\delta}(f)-f\|_X=0\}$.

\begin{theorem}\cite{recent}\label{zer}
	Let $X$ be a Banach function space and $1\in X_a$. Then the set $C_0^{\infty}([0,1])$ is dense in $X^S$.
\end{theorem}
\begin{remark}
	We observe that if $1\in X_a$, the subspace $X^S$ contains $C[0,1]$. This can be seen from the following arguments.
	$1\in X_a$ will imply  that $1\in X^S$ and hence all the constant functions are in $X^S$.
	Let $f\in C[0,1]$, consider $g(x)=f(x)-f(1)$.
	Then \[T_{\delta}g(x)-g(x)=\begin{cases}f(x+\delta)-f(x)& if\ x+\delta\in [0,1]\\ f(1)-f(x)& if\ x+\delta\not\in [0,1]\end{cases}\]\\
	For a given $\epsilon>0$, using continuity of $f$ at $x$ and at $1,$ we get a $\delta>0$ such that 
	\[
	|T_{\delta}g(x)-g(x)|=|f(x+\delta)-f(x)|<\epsilon, \textrm{ whenever  } x+\delta\in[0,1]\]\[
	\textrm{ and  } |T_{\delta}g(x)-g(x)|=|f(1)-f(x)|<\epsilon,
	\textrm{whenever }  x+\delta\not\in [0,1].\]
	Therefore $|T_{\delta}g(x)-g(x)|<\epsilon$ for every $x\in [0,1],$ as $\delta\downarrow0.$
	And thus for a given $\epsilon$,  $\rho(|T_{\delta}g(x)-g(x)|)<\epsilon\rho(1)$,
	or $\|T_{\delta}g-g\|_X<\epsilon\|1\|_X$ as $\delta\downarrow 0.$
	Hence $g\in X^S$ which gives $f \in X^S$. 
\end{remark}
The following is a Korovkin-type theorem in this setting.
\begin{theorem}\cite{recent}\label{zer1}
	Let $X$ be a Banach function space such that $1\in X_a$ and $\{L_n\}_{n\in \mathbb{N}}$ be a sequence of positive bounded linear operators in $X^S$ satisfying the condition 
	\[
	\lim\limits_{n\rightarrow\infty}L_n(g)=g\ in\ C[0,1]\ for\ all\ g\in\{1,t,t^2\}
	\]
	Then the relation $\lim\limits_{n\rightarrow\infty}L_n(f)=f$ is true in $X$ for all $f\in X^S$ if and only if $sup_{n}\|L_n\|_{B(X^S)}=c<+\infty$.
\end{theorem}

Let $X$ be a Banach function space. We denote by $X_{2\pi}^S$ the Banach space of all measurable functions $f$ such that $f\in X^S[-\pi,\pi]$ and satisfies $f(x+2\pi)=f(x)$ for a.e. $x\in\mathbb{R}$. The trigonometric analogue (the test functions being $1,\sin x, \cos x$) of the above theorem is the following.
\begin{theorem}\cite{recent}\label{zer2}
	Let $X$ be a Banach function space such that $1\in X_a$, and let $\{L_n\}_{n\in \mathbb{N}}$ be a sequence of positive bounded linear operators in $X_{2\pi}^S$ satisfying the condition 
	\[
	\lim\limits_{n\rightarrow\infty}L_n(g)=g\ in\ C_{2\pi}(\mathbb{R})\  for\ all\ g\in\{1,\sin x, \cos x\} 
	\]
	Then the relation $\lim\limits_{n\rightarrow\infty}L_n(f)=f $ holds  for all $f\in X_{2\pi}^S$
	if and only if $\sup\limits_n \|L_n\|_{B(X_{2\pi}^S)}=c<\infty$.
\end{theorem}

In this article, we present three significant results within the framework of Banach function spaces. Firstly, we extend the theorem to encompass operators that are non-positive. We provide an example of such a non-positive operator where the theorem remains applicable. Furthermore, we derive an operator version of the result using Dumitru Popa's operator version of the Korovkin Theorem obtained in \cite{popa}, along with illustrative examples. Additionally, akin to the quantitative results obtained by Shisha and Mond, De Vore, and others \cite{devoretext, nishishiraho, nishishiraho1, devore}, we establish quantitative results for the Korovkin-type theorem in \cite{recent}. Our results will be helpful to estimate the rate of convergence in the approximation results obtained in \cite{recent}. The quantitative forms in the case of the positive Kantorovich polynomials on Lebesgue space, Grand Lebesgue space etc., are also obtained.  It is noteworthy that in the Lebesgue space, Swetits and Wood \cite{swetits} obtained sharp quantitative estimates in terms of the second-order modulus of smoothness. However, our quantitative estimates extend to general Banach function spaces, covering a dense class. For instance, in the weighted Lebesgue space, we achieve significantly improved estimates by selecting suitable weights (See Remark $2.21$). We numerically demonstrate this enhancement in the last section.\par
The article is organized as follows: In the following section, we explore the main results. Initially, in subsection one, we illustrate that the assumption of positivity is not necessary, supported by an illustrative example. Subsequently, we derive an operator version of the result outlined in \cite{recent}. Progressing further, we present the quantitative Korovkin-type theorem within this context. Examples and numerical illustrations are provided in the final subsection.

\section{Main Results}
In this section, three major results pertaining to the Korovkin-type theorems Theorem \ref{zer1} and \ref{zer2} are proved. We establish that the positivity assumption in the theorem is not a necessary requirement and obtain a modified Korovkin-type theorem. Subsequently, we prove an operator version similar to Theorem \ref{popa}. Finally, we obtain a quantitative form of these theorems. We apply our results to appropriate examples and provide numerical illustrations for specific cases.
\subsection{ A Non Positive Version}
In this subsection, we establish that in the Korovkin-type theorem obtained in \cite{recent}, the positivity of the sequence of operators is not necessary. We provide an example to illustrate this. To begin, we recall the following theorems from \cite{wulbert} and \cite{nonpositive} respectively.

\begin{theorem}\cite{wulbert}\label{wul}
	Let $\{L_n\}$ be a sequence of norm one operators defined on $C[0, 1]$ ($C_{2\pi}(\mathbb{R})$, respectively). Then $L_n(f)$ converges to $f$ uniformly for all $f$ in $C[0, 1]$ ($C_{2\pi}(\mathbb{R})$, respectively) if and only if $L_n(p)$ converges to $p$ for the three functions $1$, $x$, and $x^2$ ($1$, $\cos x$ and $\sin x$ respectively). 
\end{theorem}
\begin{remark}
	Note that Theorem \ref{wul} is applicable for sequence of operators $\{L_n\}$ whose operator norm less than or equal to one. 
\end{remark}
\begin{theorem}\cite{nonpositive}\label{vin}
	Let $\{L_n\}$ be a sequence of operators defined on $C[0, 1]$ ($C_{2\pi}(\mathbb{R})$, respectively) such that $\lim\limits_n\|L_n\|=\alpha$ where $\alpha\geq 1$. Then $L_n(f)$ converges to $f$ uniformly for all $f$ in $C[0, 1]$ ($C_{2\pi}(\mathbb{R})$, respectively) if and only if $L_n(p)$ converges to $p$ for the three functions $1$, $x$, and $x^2$ ($1$, $\cos x$ and $\sin x$ respectively). 
\end{theorem}

We obtain the following Theorem:
\begin{theorem}\label{nonpositive}
	Let $X$ be a Banach function space such that $1\in X_a$ and $\{L_n\}_{n\in \mathbb{N}}$ be a sequence of bounded linear operators in $X^S$ satisfying the conditions: 
	\begin{enumerate}
		\item $\lim\limits_{n\rightarrow\infty}L_n(g)=g\ in\ C[0,1]\ for\ all\ g\in\{1,t,t^2\}$\\
		\item $\sup\limits_{n}\|L_n\|_{B(C[0,1])}<\infty$.
	\end{enumerate}
	Then the relation $\lim\limits_{n\rightarrow\infty}L_n(f)=f$ is true in $X$ for all $f\in X^S$ if and only if $\sup\limits_{n}\|L_n\|_{B(X^S)}< \infty$.
\end{theorem}

\begin{proof}
	Suppose that $\lim\limits_{n\to\infty} L_n f = f$ for all $f \in X^S$ in $X$. Then, by the uniform boundedness principle, we have $\sup_n \| L_n \|_{\mathcal{B}(X^S)} = c < +\infty$. \par 
	Conversely, suppose that there exist an $f\in C[0,1]$ such that $L_n(f)$ does not converge to $f$ as $n\to \infty$, then there exists a subsequence $\{L_{n_k}(f)\}$ of $\{L_{n}(f)\}$ and an $\epsilon >0$ such that for all $k\in \mathbb{N}$, 
	\begin{equation}
		\|L_{n_k}(f)-f\|_\infty>\epsilon
	\end{equation}
	Since $\sup\limits_k \| L_{n_k} \|_{\mathcal{B}(C[0,1])} = c < +\infty$, there exists a subsequence say $\{L_{n_l}\}$ of $\{L_{n_k}\}$ such that $\lim\limits_{l\to\infty} \| L_{n_l} \|_{\mathcal{B}(C[0,1])}=c$. From Theorem \ref{wul} and \ref{vin}, since $\lim\limits_{l\to\infty} L_{n_l}(f)=f$ for $f\in \{1,x,x^2\}$, we have 
	\begin{equation*}
		\lim\limits_{l\to\infty} L_{n_l}(f)=f \ for\ all\ f\in C[0,1] 
	\end{equation*}
	But this is a contradiction to $(2)$. Thus $\lim\limits_{n\to\infty} L_{n}(f)=f$ for all $f\in C[0,1]$.\par
	Suppose that $f \in X^S$ and $\varepsilon > 0$ be given. Then from Theorem \ref{zer} it follows that there exists $g \in C[0, 1]$ such that 
	\[
	\| f - g \|_X < \varepsilon. 
	\]
	Then $\lim\limits_{n\to\infty} L_n g = g$ in $C[0, 1]$. Hence, there exists the number $n_{\varepsilon}$ such that for all $n > n_{\varepsilon}$ 
	\[
	\| L_n g - g \|_{\infty} < \varepsilon. 
	\]
	We have, $\| g \|_X \leq c_0 \| g \|_{\infty}$, where $c_0 = \| 1 \|_X$. Now
	\begin{align*}
		\| L_n f - f \|_X &\leq \| L_n f - L_n g \|_X + \| L_n g - g \|_X + \| f - g \|_X \\
		&\leq (c + 1) \| f - g \|_X + c_0 \| L_n g - g \|_{\infty} < (c + 1)\varepsilon + c_0\varepsilon = c_1\varepsilon.
	\end{align*}
	Thus, $\lim\limits_{n\to\infty} L_n f = f$ for all $f\in X^S$. Hence the proof.
\end{proof}
The following is the trigonometric analogue of this result.
\begin{theorem}
	Let $X$ be a Banach function space such that $1\in X_a$ and $\{L_n\}_{n\in \mathbb{N}}$ be a sequence of bounded linear operators in $X_{2\pi}^S$ satisfying the conditions: 
	\begin{enumerate}
		\item $\lim\limits_{n\rightarrow\infty}L_n(g)=g\ in\ C_{2\pi}(\mathbb{R})\ for\ all\ g\in\{1,\sin(.),\cos(.)\}$\\
		\item $\sup\limits_{n}\|L_n\|_{B(C_{2\pi}(\mathbb{R}))}<\infty$
	\end{enumerate}
	Then the relation $\lim\limits_{n\rightarrow\infty}L_n(f)=f$ is true in $X$ for all $f\in X_{2\pi}^S$ if and only if $\sup\limits_{n}\|L_n\|< \infty$.
\end{theorem}
The proof follows by the same techniques as the previous theorem.
\begin{example}
	We give an example of a non positive operator for which Theorem \ref{nonpositive} is applicable. Firstly, we recall the definition of Lagrange Interpolation Operator.
	\begin{definition}
		Let $\{x_1,x_2 \ldots x_n\}$ be a set of $n$ nodal/grid points in $[-1,1]$ and suppose $f\in C[-1,1]$, then the Lagrange interpolation operator gives a unique polynomial of degree $(n-1)$ assuming the values $\{f(x_1), f(x_2),\ldots,f(x_n)\}$ at the nodal/grid points $x_1,x_2,\ldots x_n$ respectively which is given by,
		\[
		L_n(f)(x)=\sum\limits_{k=1}^{n}f(x_k)P_k(x)   
		\]
		$x\in[-1,1]$. Here 
		\[
		P_k(x)=\frac{\zeta(x)}{\zeta\prime(x_k)(x-x_k)}
		\]
		$k=1,\ldots,n$ and the polynomials $\omega(x)$ is defined by 
		\[
		\zeta(x)=c(x-x_1)(x-x_2)\ldots (x-x_n)
		\]
		where $c$ is an arbitrary non-zero constant.
	\end{definition} 
	\begin{remark}
		If we choose the nodal points to be
		\[
		x_k=\cos \theta_k^{(n)}, where\ \theta_k^{(n)}=\frac{(2k-1)\pi}{2n}, \ k=1,2,\dots n
		\]
		(the Chebychev nodes of first kind), then 
		\[
		\zeta(x)=T_n(x)=\cos(n\arccos x)=\cos n\theta,\ \cos\theta=x
		\]
		$\zeta^\prime(\theta)=\frac{n\sin n\theta}{\sin\theta}$, so that 
		\[
		P_k(\theta)=(-1)^{k+1}\frac{\cos n\theta\sin \theta_k^{(n)}}{n(\cos\theta-\cos \theta_k^{(n)})}
		\]
		$k=1,\ldots n$ and
		\[
		L_n(f)(\theta)=\sum\limits_{k=1}^nf(\cos\theta_k^{(n)})(-1)^{k+1}\frac{\cos n\theta\sin \theta_k^{(n)}}{n(\cos\theta-\cos \theta_k^{(n)})}
		\]
	\end{remark}
	Consider the operator $H_{n,\delta}:L^1[0,\pi] \rightarrow L^1[0,\pi]$ defined as   
	\[
	H_{n,\delta}(f)(x)=\frac{1}{2\pi}\sum_{l=-\infty}^{\infty}\int\limits_{l\pi}^{(l+1)\pi}\frac{1}{2}\sum\limits_{k=1}^n\widehat{f\chi_{[0,\pi]}*\phi_\delta}(\theta_k^{(n)}+l\pi)(P_k^l(\theta+\frac{\pi}{2n})+P_k^l(\theta-\frac{\pi}{2n}))e^{ix\theta}d\theta
	\]
	where $\hat{f}$ denotes the Fourier transform of $f$ on $\mathbb{R}$ and $\{\phi_\delta\}$ is an approximate identity. This operator has been constructed in \cite{nonpositive}. It is seen that this operator is not positive for specific choice of $\phi_\delta$ page $29$ Remark $6.4$ \cite{nonpositive}. We also have 
	\[
	\|H_{n,\delta}(f)\|_1\leq C_\delta\|f\|_1
	\]
	where $C_\delta$ is a constant which does not depend on $n$ or $f$. Moreover, $\sup\limits_n\|H_{n,\delta}\|_{B(C[0,\pi])}<\infty$. For more details refer \cite{nonpositive}. Recall that $(L^1[0,\pi],\|.\|_1)$ is a Banach function space. In this example, we consider the test functions $\{1,\cos(.),\cos^2(.)\}$. Hence it suffices to prove the convergence on $\{1, \cos(.), \cos^2(.)\}$ by Theorem \ref{nonpositive}. We have established this convergence in \cite{nonpositive}(page $35$). Hence applying Thoerem \ref{nonpositive}, we have $\lim\limits_{n\to\infty, \delta\to 0}H_{n,\delta}(f)=f$ for all $f\in (L^1[0,\pi])^S$.
\end{example}
\subsection{An Operator Version}
In this subsection, we demonstrate that for $X$, a real Banach function space, Theorem \ref{zer1} can be further generalized to an operator version of the theorem. This result aligns with an operator version of the Korovkin theorem obtained by Dumitru Popa \cite{popa}. Popa extended the theorem by substituting the convergence of the sequence of positive linear operators to the identity operator with convergence to an arbitrary operator, and he proved sufficient conditions under which a Korovkin-type theorem is obtained. We state the result below. \par 
Let $e_i(t)=t^i$ for $t\in[a,b]$ and $i=0,1,2$. We state the result below.
\begin{theorem}\cite{popa}\label{popa}
	Let $T$ be a compact Hausdorff space, $V_n:C[a,b]\rightarrow C(T)$ a sequence of positive linear operators and $A:C[a,b]\rightarrow C(T)$ a linear operator such that $A(e_0)A(e_2)=[A(e_1)]^2$ and $A(e_0)(t)>0$ for every $ t\in T$. If $\lim\limits_{n\rightarrow\infty}V_n(e_i)=A(e_i)$ for $i=0,1,2$ all uniformly on $T$, then for every $f\in C[a,b]$, $\lim\limits_{n\rightarrow\infty}V_n(f)=A(f)$ uniformly on $T$.
\end{theorem} 
In \cite{popa2}, the author derived the trigonometric analogue of this result, accompanied by several illustrative examples. Moreover, in \cite{vin}, we obtained a quantitative form of this result, thereby obtained the trigonometric analogue as a consequence. Subsequently non linear operator version and its quantitative forms of Theorem \ref{popa} were obtained in \cite{gal} and \cite{gal2} respectively. Recently, in \cite{dorai}, the author has obtained a characterization for the condition $A(e_0)A(e_2)=[A(e_1)]^2$ and $A(e_0)(t)>0$ for every $ t\in T$. He proved that in this case the operator $A$ will be the weighted composition operator given by $A(f)=A(e_0)f\circ (\frac{A(e_1)}{A(e_0)})$ for every $f\in C[a,b]$.\par
We establish that Theorems \ref{zer1} and \ref{zer2} can be further extended to an operator version of the theorem. The following results are obtained:
\begin{theorem}\label{oper1}
	Let $X$ be a real Banach function space such that $1\in X_a$ and $\{L_n\}_{n\in \mathbb{N}}$ be a sequence of positive bounded linear operators on $X^S$ and $A$ be a positive bounded linear operator in $X^S$ satisfying the conditions:
	\begin{enumerate}
		\item{$A(e_0)A(e_2)=[A(e_1)]^2$ and $A(e_0)(t)>0$ for every $ t\in [0,1]$}
		\item{	
			$\lim\limits_{n\rightarrow\infty}L_n(g)=A(g)\ in\ C[0,1]\ for\ all\ g\in\{e_0,e_1,e_2\}$}
	\end{enumerate}
	Then, $\lim\limits_{n\rightarrow\infty}L_n(f)=A(f)$ in $X$ for all $f\in X^S$ if and only if $\sup\limits_{n}\|L_n\|< \infty$.
\end{theorem}
\begin{proof}
	Let $\lim\limits_{n\to\infty} L_n f = A f$ for all $f \in X^S$ in $X$. Is is also given that $A$ is a bounded linear operator. Then, by the uniform boundedness principle we must have, $\sup_n \| L_n \|_{\mathcal{B}(X^S)} < +\infty$.
	
	To prove the sufficient condition, consider a function $f \in X^S$ and an $\varepsilon > 0$. Then from Theorem \ref{zer} it follows that there exists $g \in C[0, 1]$ such that 
	\[
	\| f - g \|_X < \varepsilon. \tag{6}
	\]
	By Theorem \ref{popa}, we have $\lim\limits_{n\to\infty} L_n g = Ag$ in $C[0, 1]$. Consequently, there exists the number $n_{\varepsilon}$ such that for all $n > n_{\varepsilon}$ 
	\[
	\| L_n g -A g \|_{\infty} < \varepsilon. \tag{7}
	\]
	Since, $\| g \|_X \leq c_0 \| g \|_{\infty}$, where $c_0 = \| 1 \|_X$, and by $(6)$ and $(7)$,
	\begin{align*}
		\| L_n f - Af \|_X &\leq \| L_n f - L_n g \|_X + \| L_n g - A g \|_X + \|A f - A g \|_X \\
		&\leq (c + \tilde{c}) \| f - g \|_X + c_0 \| L_n g - A g \|_{\infty} < (c + \tilde{c})\varepsilon + c_0\varepsilon = c_1\varepsilon.
	\end{align*}
	Hence, $\lim\limits_{n\to\infty} L_n f = A f$ for all $f\in X^S$.
\end{proof}
\begin{example}
	The Kantorovich polynomial is given by
	\begin{equation}
		K_n(f)(x)=(n+1)\sum_{k=0}^{n}\binom{n}{k} x^k (1-x)^{n-k}\int_{\frac{k}{n+1}}^{\frac{k+1}{n+1}}f(t) dt
	\end{equation}
	$n=1,2,\ldots$, where $x\in[0,1]$ and $f\in L^1[0,1]$.\\
	Note that $K_n(1)=1$ for $n=1,2,...$\par
	Consider the modified Kantorovich operators defined by,
	\begin{equation}
		K_n^{\alpha,\beta}(f)(x)=(n+1)\sum_{k=0}^{n}\binom{n}{k} x^k (1-x)^{n-k}\int_{\frac{k}{n+1}}^{\frac{k+1}{n+1}}f(\alpha t+ \beta) dt
	\end{equation}
	$n=1,2,\ldots$, where $x\in[0,1]$, $\alpha, \beta$ are constants and $f\in L^1[0,1]$.\par
	We apply Theorem \ref{oper1} to these operators. Note that the uniform boundedness of the sequence of operators $\{K_n\}$ is proved in \cite{recent} and hence the operators $\{K_n^{\alpha,\beta}\}$ is uniformly bounded. Let $A(f)(x)=f(\alpha x +\beta )$. Then $A$ is a positive bounded linear operator in $L^1[0,1]$. Moreover $A$ satisfies the desired conditions.  Hence by Theorem \ref{oper1}, $\lim\limits_{n\to\infty}K_n^{\alpha,\beta}(f)=A(f)$ for all $f\in L^1[0,1]^S$.
\end{example}
Now we state the trigonometric analogue of this theorem. We denote $h_0(x)=1$, $h_1(x)=\cos x$ and $h_2=\sin x$ for all $x\in \mathbb{R}$.
\begin{theorem}\label{oper2}
	Let $X$ be a real Banach function space such that $1\in X_a$ and $\{L_n\}_{n\in \mathbb{N}}$ be a sequence of positive bounded linear operators on $X_{2\pi}^S$ and $A$ be a positive bounded linear operator in $X_{2\pi}^S$ satisfying the conditions:
	\begin{enumerate}
		\item{$A(h_0)^2=A(h_1)^2+A(h_2)^2$ and $A(h_0)(t)>0$ for every $ t\in [-\pi,\pi]$}
		\item{	
			$\lim\limits_{n\rightarrow\infty}L_n(g)=A(g)\ in\ C[0,1]\ for\ all\ g\in\{h_0,h_1,h_2\}$}
	\end{enumerate}
	Then $\lim\limits_{n\rightarrow\infty}L_n(f)=A(f)$ in $X$ for all $f\in X_{2\pi}^S$ if and only if $\sup\limits_{n}\|L_n\|< \infty$.
\end{theorem}
The proof is similar to the proof of Theorem \ref{oper1}.
\begin{example}
	The following example is an extension of the sequence of operators considered in Corollary $2$, \cite{popa2}. We extend the definition from $C_{2\pi}[-\pi,\pi]$ to $L^1[-\pi,\pi]$.
	Consider $L_n: L^1[-\pi,\pi]\to L^1[-\pi,\pi]$ defined by (See \cite{popa2}),
	\[
	L_n(f)(x)=\sum\limits_{k=0}^{n}\binom{n}{k}\sin^{2k} x\cos^{2{n-k}}x \int\limits_{0}^{1}f(\frac{t+k}{n+1}+x)dt
	\]
	Then $\{L_n\}$ is uniformly bounded in the operator norm. Let $A(f)(x)= f(x+\sin^2 x)$. Then we have $A(h_0)^2=A(h_1)^2+A(h_2)^2$. Moreover we also have $\lim\limits_{n\to\infty} L_n(h_i)=A(h_i)$ for $i=0,1,2$ in $C_{2\pi}(\mathbb{R})$ (Corollary $2$, \cite{popa2}). If $A$ is bounded then we are done. If not, let $w(x)$ be any weight function such that $A(f)(x)=w(x)f(x+\sin^2 x)$ is bounded in $L^1[0,1]$. For instance, let $w(x)=1+\sin 2x$. Then $\|A(f)\|_1\leq c \|f\|_1$ for some constant $c$. Thus we can modify our sequence of operators to $\{w L_n\}$ and it satisfies the conditions of Theorem \ref{oper2}. 
\end{example}
\subsection{Quantitative Results}
In this section, we prove the quantitative forms of the Korovkin-type theorems obtained in \cite{recent}.
We obtain results analogous to the results in \cite{shisha, devore}. To begin with we recall the quantitative result obtained by Shisha and Mond below:
\begin{theorem}{\cite{shisha}}\label{shih1}
	Let $\{L_n\}_{n\in\mathbb{N}}$ be a sequence of positive linear operators with the same domain $D$ which contains the restrictions of $1,t,t^2$
	to $[a,b]\subset \mathbb{R}$. For $n=1,2,\ldots$, suppose $L_n(1)$ is bounded. Let $f\in D$ be continuous in $[a,b]$, with modulus of continuity $\omega$. With $\|.\|$ denotes the  sup-norm, we have for  $n=1,2,\ldots$,
	\[
	\|f-L_n(f)\|\leq\|f\|\|L_n(1)-1\|+ \|L_n(1)+1\|\omega(f,\mu_n)
	\textrm{ where } 
	\mu_n^2:=\|L_n((t-x)^2)(x)\|. 
	\]
	In particular, if $L_n(1)=1$, this inequality becomes 
	$\|f-L_n(f)\|\leq2\omega(f,\mu_n)$
\end{theorem}
In this result we have the convergence of the sequence of positive linear operators to the identity operator given in terms of the convergence on the test functions $\{1, x, x^2\}$ and the modulus of continuity of the function. The classical modulus of continuity is defined as follows:
\begin{definition}
	Let $f:I\to \mathbb{R}$ denote a bounded function on a real interval $I$. Then the modulus of continuity of $f$ with argument $\delta>0$ is defined as 
	\[
	\omega(f,\delta):=sup\{|f(x)-f(y)|: x,y\in I,\ |x-y|\leq\delta\}
	\]
\end{definition}
The following are few properties of the modulus of continuity:
\begin{enumerate}
	\item {If $0<\delta_1<\delta_2$, then $
		\omega(f, \delta_1)\leq \omega(f,\delta_2).$ }
	\item{Let $\delta>0, \,\lim\limits_{\delta\rightarrow 0+}\omega(f,\delta)=0$
		for all uniformly continuous and bounded functions $f$ on $I$.}
	\item{ Let $\lambda, \delta >0$, then $
		\omega(f,\lambda\delta)\leq(1+[\lambda])\omega(f,\delta)$			where $[.]$ denotes the integer part of the function.}
\end{enumerate}
The following inequality was crucial to obtain the quantitative form which follows from the properties of the modulus of continuity.:
For $x,y\in [a,b]$, $f\in C[a,b]$, $\delta>0$,
\begin{equation*}
	|f(x)-f(y)|\leq (1+(x-y)^{2}\delta^{-2})\omega(f,\delta)
\end{equation*}
\par
Furthurmore, they proved the trigonometric analogue of this result.
Such quantitative theorems were studied extensively by many mathematicians and they were successful in extending these results into several Korovkin-type settings. For instance, refer \cite{devoretext, nishishiraho, nishishiraho1}. In $1972,$ Ronald A DeVore obtained an extension to this quantitative version while exploring the optimal conditions in \cite{devore} and obtained saturation theorems in the case of positive polynomial operators.

Analogous to the  quantitative versions of Theorems \ref{kor1} and \ref{kor2} and the quantitative theorem Theorem \ref{shih1}, we prove the	quantitative versions of Theorems \ref{zer1}, \ref{zer2}. We adopt the techniques from \cite{shisha, devore} to prove this theorem. 
\begin{theorem}\label{quan1}
	Let $\{L_n\}_{n\in\mathbb{N}}$ denote a sequence of positive linear operators in $X^S$ and let $1\in X_a$. Suppose that $\sup\limits_n\|L_n(1)\|_X=c<+\infty$. Let $\omega(f,\delta)$ be the modulus of continuity of $f\in C[0,1]$ with the argument $\delta>0$. Then for every $f\in C[0,1]$,
	\[
	\|L_n(f)-f\|_X\leq\|f\|_{\infty}\|L_n(1)-1\|_X + (\|L_n(1)\|_X+1)\omega(f,\mu_n),\ n=1,2,3,\ldots, 
	\]
	where $
	\mu_n^2=\|L_n((x-.)^2)(x)\|_X.$		If $L_n(1)=1$ then we have
	\[
	\|L_n(f)-f\|_X\leq (c+1)\omega(f,\mu_n), n=1,2,3,\ldots
	\]
	Moreover, if $f$ is differentiable and $f'\in C[0,1]$ then we have 
	\[
	\|L_n(f)-f\|_X\leq\|f\|_{\infty}\|L_n(1)-1|)\|_X+\sqrt{c}\mu_n\|f'\|_{\infty} +(\sqrt{c}+1)\mu_n\omega(f',\mu_n).
	\]
	In addition, if $L_n(1)=1$, we get $
	\|L_n(f)-f\|_X\leq \sqrt{c}\mu_n\|f'\|_{\infty}+ (\sqrt{c}+1)\mu_n\omega(f',\mu_n).$
\end{theorem}
\begin{proof}
	To begin with we let $x,y\in [0,1]$ and let $f\in C[0,1]$ then 
	\begin{equation*}
		|f(x)-f(y)|\leq (1+(x-y)^2\delta^{-2})\omega(f,\delta)
	\end{equation*}
	which follows from the properties of the modulus of continuity, $\omega$.\\
	Now we write $	L_n(f)(x)-f(x)=L_n(f)(x)-L_n(f(x).1)+f(x)L_n(1)(x)-f(x)$
	Thus 
	\begin{equation*}
		|L_n(f)(x)-f(x)|\leq|L_n(f-f(x).1)(x)|+|f(x)||L_n(1)(x)-1|
	\end{equation*}
	By the positivity of the operator $L_n$,
	\begin{equation*}
		\begin{split}
			|L_n(f-f(x).1)|&\leq L_n(|f-f(x).1|)\\
			&\leq (L_n(1)+ L_n((x-.)^2)\delta^{-2})\omega(f,\delta)
		\end{split}
	\end{equation*}
	Thus 
	\begin{equation*}
		\begin{split}
			|L_n(f)(x)-f(x)|&\leq L_n(|f-f(x).1|)(x)+|f(x)||L_n(1)(x)-1|\\
			&\leq (L_n(1)(x)+ L_n((x-.)^2)(x)\delta^{-2})\omega(f,\delta)+ \|f\|_{\infty}|L_n(1)(x)-1|\\
			|L_n(f)-f| &\leq (L_n(1)+L_n((x-.)^2)\delta^{-2})\omega(f,\delta)+ \|f\|_{\infty}|L_n(1)-1|
		\end{split}
	\end{equation*}
	Recall the properties of the function norm we have already defined. This gives that,
	\begin{equation*}
		\begin{split}
			\rho(|L_n(f)-f|)&\leq \rho(L_n(|f-f(x).1|)+|f(x)|\rho(L_n(1)-1|)\\
			&\leq \rho{(L_n(1)+L_n((x-.)^2)\delta^{-2})}\omega(f,\delta)+ \|f\|_{\infty}\rho(|L_n(1)-1|)
		\end{split}
	\end{equation*}
	We choose $\delta^2=\mu_n^2=\|L_n((x-.)^2)(x)\|_X$ and we have $\rho(|f|)=\|f\|_X$. So that
	\begin{equation*}
		\begin{split}
			\|L_n(f)-f\|_X&\leq \|L_n(1)+L_n((x-.)^2)\delta^{-2})\|_X\omega(f,\delta) + \|f\|_{\infty}\|L_n(1)-1\|_X\\
			&\leq(\|L_n(1)\|_X+\|L_n((x-.)^2)\|_X\delta^{-2})\|_X)\omega(f,\delta) + \|f\|_{\infty}\|L_n(1)-1\|_X\\
			&\leq(\|L_n(1)\|_X+1)\omega(f,\mu_n) + \|f\|_{\infty}\|L_n(1)-1\|_X
		\end{split}
	\end{equation*}
	Let $f$ be also differentiable with $f'\in C[0,1]$ and $\delta>0$, then if $x,y\in [0,1]$ with $x<y$ the Mean Value Theorem gives us the existence of a $c\in (x,y)$ such that
	\begin{equation*}
		f(y)-f(x)=f'(c)(y-x)
	\end{equation*}
	Now 
	\begin{equation*}
		\begin{split}
			|f(y)-f(x)-f'(x)(y-x)|&\leq|y-x||f'(c)-f'(x)|\\
			&\leq|y-x|(1+|c-x|\delta^{-1})\omega(f',\delta)\\
			&\leq|y-x|(1+|y-x|\delta^{-1})\omega(f',\delta)\\
			&\leq(|y-x|+(y-x)^2\delta^{-1})\omega(f',\delta)
		\end{split}
	\end{equation*}
	The positivity of $L_n$ and the above inequality gives

	\begin{multline*}
		|L_n(f)(x)-f(x)|
		\leq L_n(|f-f(x).1-f'(x)(.-x)|)(x)+\\
		\|f'\|_{\infty}L_n(|.-x|) +\|f\|_{\infty}|L_n(1)(x)-1|
	\end{multline*}
	\begin{multline*}
		|L_n(f)(x)-f(x)|\leq	 (L_n(|.-x|)(x)+L_n((.-x)^2)(x)\delta^{-1})\omega(f',\delta)\\
		+\|f'\|_{\infty}L_n(|.-x|)+\|f\|_{\infty}|L_n(1)(x)-1|
	\end{multline*}
	By the Cauchy-Schwartz inequality for positive operators,
	\begin{equation*}
		L_n(|.-x|)^2\leq L_n((.-x)^2).L_n(1)
	\end{equation*}
	We choose $\delta=\mu_n$ then,
	\begin{multline*}
		|L_n(f)-f| \leq (\sqrt{L_n((.-x)^2)}\sqrt{L_n(1)}+L_n((.-x)^2)\mu_n^{-1}) \omega(f',\mu_n)+ \\
		\|f'\|_{\infty}(\sqrt{L_n((.-x)^2)}\sqrt{L_n(1)})+\|f\|_{\infty}|L_n(1)-1|
	\end{multline*}
	We can see that the following property of function norm, $\rho$ is true: Let $1\leq p, q<+\infty$ satisfies that $\frac{1}{p}+\frac{1}{q}=1$, then for $f,g\in C[0,1]$,
	\begin{equation*}
		\rho(|fg|)\leq\rho(|f|^p)^{\frac{1}{p}}\rho(|g|^q)^{\frac{1}{q}}
	\end{equation*}
	By this property of the function norm $\rho$ and by the Cauchy-Schwartz Inequality of the positive operators we have,\\
	\begin{multline*}
		\rho(|L_n(f)-f|)\leq (\rho(\sqrt{L_n((.-x)^2)}(\sqrt{L_n(1)}) + \rho(L_n((.-x)^2)\mu_n^{-1}))\omega(f',\mu_n) + \\ \|f'\|_{\infty} \rho(\sqrt{L_n((.-x)^2)}(\sqrt{L_n(1)})+\|f\|_{\infty}\rho(|L_n(1)-1|)
	\end{multline*}
	\begin{multline*}
		\rho(|L_n(f)-f|)\leq(\rho(L_n((.-x)^2))^{\frac{1}{2}}\rho(L_n(1))^{\frac{1}{2}} + \rho(L_n((.-x)^2)\mu_n^{-1})\omega(f',\mu_n)+\\
		\|f'\|_{\infty}(\rho(L_n((.-x)^2))^{\frac{1}{2}}\rho(L_n(1))^{\frac{1}{2}}+	\|f\|_{\infty}\rho(|L_n(1)-1|)
	\end{multline*}
	\begin{multline*}
		\|L_n(f)-f\|_X \leq(\|L_n((.-x)^2)\|_X^{\frac{1}{2}}\|L_n(1)\|_X^{\frac{1}{2}}+\|L_n((.-x)^2)\|_X\mu_n^{-1})\omega(f',\mu_n)+\\
		\|f'\|_{\infty}\|L_n((.-x)^2)\|_X^{\frac{1}{2}}\|L_n(1)\|_X^{\frac{1}{2}}+	\|f\|_{\infty}\|L_n(1)-1\|_X
	\end{multline*}
	Hence we obtain the desired inequality.
\end{proof}

The main advantage of this theorem over the one by Shisha and Mond \cite{shisha} is that it provides inequalities in the Banach space norm, which is generally weaker than the sup-norm since, $\|g\|_X\leq c_0\|g\|_{\infty}$, where $c_0=\|1\|_X$. Convergence with respect to the sup-norm guarantees the convergence with respect to the Banach function space norm. Moreover, to obtain the quantitative results, the assumption of the uniform boundedness of the sequence of operators is not necessary.  
The trigonometric analogue of is as follows:
\begin{theorem}\label{quan2}
	Let $X$ be a Banach function space such that $1\in X_a$ and let $\{L_n\}_{n\in\mathbb{N}}$ be a sequence of positive linear operators in $X_{2\pi}^S$. Suppose that $\sup\limits_n\|L_n(1)\|_X=c<+\infty$ and let $\omega(f,\delta)$ denote the modulus of continuity of $f$ with argument $\delta>0$. Then for $f\in C_{2\pi}[-\pi,\pi],$ 
	\[
	\|L_n(f)-f\|_X\leq\|f\|_{\infty}\|L_n(1)-1\|_X+ (\|L_n(1)\|_X+1)\omega(f,\mu_n),\ n=1,2,3,\ldots
	\]
	where $\mu_n^2=\pi^2\|L_n(sin^2\frac{(x-.)}{2})(x)\|_X$. If $L_n(1)=1$ then,
	\[
	\|L_n(f)-f\|_X\leq(c+1)\omega(f,\mu_n)
	\]
	Moreover if $f$ is also differentiable and $f'\in C_{2\pi}[-\pi,\pi]$ then we have 
	\[
	\|L_n(f)-f\|_X\leq\|f\|_{\infty}\|L_n(1)-1|)\|_X+\sqrt{c}\mu_n\|f'\|_{\infty} +(\sqrt{c}+1)\mu_n\omega(f',\mu_n).
	\]
	In addition, if $L_n(1)=1$, we get $
	\|L_n(f)-f\|_X\leq \sqrt{c}\mu_n\|f'\|_{\infty}+ (\sqrt{c}+1)\mu_n\omega(f',\mu_n).$
\end{theorem}
\begin{proof}
	Let $x,y \in \mathbb{R}$, $f\in C_{2\pi}([-\pi, \pi])$. We modify the previous inequality as\\
	\begin{equation*}
		\begin{split}
			|f(x)-f(y)|&\leq (1+(x-y)^2\delta^{-2})\omega(f,\delta)\\
			&\leq (1+\pi^2sin^2\frac{(x-y)}{2}\delta^{-2})\omega(f,\delta)
		\end{split}
	\end{equation*}
	Now proceeding in a similar manner as the previous theorem we obtain the desired inequality.
\end{proof}

\subsection{Examples and Numerical Illustrations}
We recall the Kantorovich polynomials defined in Section 2.2, Example 1. An analogue of Korovkin-type theorems for these polynomial operators is derived, both in the cases of rearrangement-invariant and general non-rearrangement invariant Banach function spaces, as presented in \cite{recent}. This result was obtained for Lebesgue spaces, Grand Lebesgue spaces, Morrey-type spaces, and others. The authors proved the boundedness of the Hardy littlewood maximal operator defined by
\begin{equation*}
	Mf(x):=\sup\limits_{x\in I}\frac{1}{I}\int_{I}|f(t)|dt
\end{equation*}
and hence the boundedness of the operator norm $\|K_n\|_{B(X)}$ in each of the above examples to illustrate the theorem. In our case however we do not treat these cases separately. Nevertheless, we obtain the quantitative forms in the case of the Kantorovich polynomials on these spaces. We see how these estimates differ in each case. 
\begin{example}{\textbf{Lebesgue Space}}
	
	Let $1\leq p<+\infty$ and $X=L^p(0,1)$  with the measure space $([0,1], B , \mu)$ where $\mu$ here denotes the Lebesgue measure on $[0,1]$, $B$ is the Borel $\sigma-$algebra on $[0,1]$. The function norm we consider here is the integral function norm defined by:
	\begin{equation*}
		\rho(f)=(\int_{0}^{1}f(x)^pdx)^{\frac{1}{p}} 
	\end{equation*}
	where $f$ is non-negative. Then $\rho$ satisfies all the axioms of the function norm (See Definition $1$). Clearly axiom 1 holds since the R.H.S is the $L_p-$norm of the function $f$. 2 , 3 and 4 holds by the properties of integral and the Dominated convergence theorem.\par
	For $p=1$, the axiom 5 holds. Suppose $p>1$ we have $m([0,1])=1<\infty$. Then we can find a $q$ satisfying $\frac{1}{p}+\frac{1}{q}=1$. Holder's inequality gives
	\begin{equation*}
		\int_{[0,1]}fd\mu\leq\|f\|_p.\|1\|_q=\|f\|p=\rho(f)
	\end{equation*}
	so that axiom 5 also holds.
	$\|f\|_{L_p([0,1])}=\rho(|f|)=\|f\|_p$, the $L_p$ norm.\par
	Here $X^S$ is the closure in $X$ of the subspace $\{f:\lim\limits_{\delta\rightarrow 0}\|T_{\delta}f-f\|_X=0\}$.
	Now we apply the above quantitative theorem provided all the assumptions of the Theorem \ref{quan1} hold.\par
	For $f\in C[0,1]$ we have \\
	\begin{equation*}
		\|L_n(f)-f\|_p\leq \|f\|_{\infty}\|L_n(1)-1\|_p+ (\|L_n(1)\|_p+1)\omega(f,\mu_n),\ n=1,2,...
	\end{equation*}
	where 
	\begin{equation*}
		\mu_n^2=\|L_n((x-.)^2)(x)\|_p
	\end{equation*}
	And when $L_n(1)=1$,
	\begin{equation*}
		\|L_n(f)-f\|_p\leq (c+1)\omega(f,\mu_n),\ n=1,2,3...
	\end{equation*}
	\textbf{ Case $p=1$:}\par
	Let $p=1$ so that $X=L^1[0,1]$ and let $\{K_n\}$ be the sequence of positive Kantorovich polynomials defined on $X^S$ (Section $2.2$, Example $1$). Let $1\in X_a$. Here we have for each $n=1,2,..$, $K_n(1)=1$ which gives $c=1$. Then by the $Theorem\ 9$, for $f\in C[0,1]$ we have 
	\begin{equation*}
		\|K_n(f)-f\|_1\leq 2\omega(f,\mu_n),\ n=1,2,...
	\end{equation*}
	where
	\begin{equation*}
		\mu_n^2=\|K_n((x-.)^2)(x)\|_1
	\end{equation*}
	Let $x\in [0,1]$, we evaluate $\|K_n((x-.)^2)(x)\|_1$.\\
	Note that $K_n[(x-y)^2](x)=x^2-2xK_n(y)(x)+K_n(y^2)(x)$
	\begin{equation*}
		K_n(y)(x)=\frac{2nx+1}{2(n+1)},
		K_n(y^2)(x)=\frac{3n(n-1)x^2+6nx+1}{3(n+1)^2}
	\end{equation*}
	We have 
	\[
	K_n(y-x)(x)=\frac{1-2x}{2(n+1)}, K_n((x-y)^2)(x)=\frac{1+3(n-1)x-3(n-1)x^2}{3(n+1)^2}
	\]
	After computation and simplification we get,
	\begin{equation*}
		\mu_n^2=\|K_n[(x-y)^2](x)\|_1=\frac{1}{6(n+1)}
	\end{equation*}
	Thus we obtain for $f\in C[0,1]$
	\begin{equation*}
		\|K_n(f)-f\|_1\leq 2 \omega(f,\frac{1}{\sqrt{6(n+1)}}),\ n=1,2,...
	\end{equation*}
\end{example}

For $p=2$, applying Theorem \ref{quan1}, we have 
\[
\|K_n(f)-f\|_2\leq 2\omega(f,\mu_n)
\]
where we can compute and observe that $\mu_n^2=\mathcal{O}(\frac{1}{n+1})$.
\begin{remark}
	It is important to note that in the Lebesgue space, sharp estimates in terms of second modulus of smoothness were obtained by Swetits and Wood \cite{swetits}. The following is a quantitative Korovkin-type theorem in the $L^p$-space in terms of the second modulus of smoothness obtained by them \cite{swetits}. The second order modulus of smoothness in $L^p[a,b]$ for $1\leq p< \infty$ for argument $h>0$ is defined as 
	\[
	\omega_{2,p}(f,h)=\sup\limits_{0<t\leq h}\|f(.+t)-2f(.)+f(.-t)\|_{L^p(I_{2t})}
	\]
	where $L^p(I_{2t})$ indicates the $L^p$ norm is taken over $[a + t, b - t]$. We state the theorem below.
	\begin{theorem}\cite{swetits}\label{swetits}
		Let $\{L_n\}$ be a uniformly bounded sequence of positive linear operators from $L^{p}[a, b]$ into $L^{p}[c, d]$, where $1 \leq p < \infty$, $a \leq c < d \leq b$, and assume $\mu_{np}\to 0$ ($n\to\infty$). Then for $f \in L^{p}[a, b]$ and $n$ sufficiently large,
		\[
		\|L_n f - f\|_{L^{p}[c, d]} \leq C_{p}(\mu_{np}^2\|f\|p+\omega_{2,p}(f,\mu_{np}))
		\]
		where $C_{p} > 0$ is independent of $f$ and $n$ and $\mu_{np}=(\sup\{\|L_n(e_0)-e_0\|_{L_p[c,d]},\|L_n((t-x),x)\|_{L_p[c,d]},\|L_n((t-x)^2,x)\|_{L_p[c,d]}^{\frac{2p}{2p+1}}\})^\frac{1}{2}$
	\end{theorem} 
\end{remark}
In the following part, we give a comparison of the rate of  convergence in Theorem \ref{quan1} and \ref{swetits} numerically. We compute the modulus of continuity and the second order modulus of smoothness for $p=1$ and observe the rate of convergence for different functions on $C[0,1]$.
Let $p=1$. We have, 
\[
\mu_{n,1}^2=\max\{\frac{1}{4(n+1)}, (\frac{1}{6(n+1)})^\frac{2}{3}\}
\]
Hence we have 
\[
\|K_n(f)-f\|_1\leq C_1(\mu_{n,p}^2\|f\|_1+\omega_{2,1}(f,\mu_{n,1}))
\]

\begin{enumerate}            
	\item 
	\[
	f(x)=x^2
	\]
	\begin{center}
		\begin{tabular}{ |p{1cm}|p{2.7cm}|p{2.8cm}|p{2.7cm}|p{2.8cm}|  }
			\hline
			&&&& \\
			n& $\mu_n(p=1)$&$\omega(f,\mu_n)
			(p=1)$&$\mu_{n1}$&$\omega_{2,1}(f,\mu_{n1})$\\
			\hline
			10& 0.123091490979&0.22933043153&0.2474488015&0.1649627755\\
			\hline
			100&0.04062222318&0.07847687527&0.11817051419&0.07877928433\\
			\hline
			200&0.02879561418&0.0552704857&0.093947284323&0.0626303559\\
			\hline
			300&0.023531040266&0.045515986456&0.082115930412&0.05474306209\\
			\hline
			500&0.0182391885&0.0357113870627&0.06928996488&0.046192744776\\
			\hline
			1000&0.012903494&0.023879735591&0.05501378899&0.036675211880\\
			\hline
		\end{tabular}
	\end{center}
	\item 
	\[
	f(x) = \begin{cases} 
		2x & \text{if } 0 \leq x \leq 0.5 \\
		-2x + 2 & \text{if } 0.5 < x \leq 1 
	\end{cases}
	\]
	\begin{center}
		\begin{tabular}{ |p{1cm}|p{2.7cm}|p{2.8cm}|p{2.7cm}|p{2.8cm}|  }
			\hline
			&&&& \\
			n& $\mu_n(p=1)$&$\omega(f,\mu_n)(p=1)$&$\mu_{n1}$&$\omega_{2,1}(f,\mu_{n1})$\\
			\hline
			10& 0.123091490979&0.24424424424&0.2474488015&0.08206989775\\
			\hline
			100&0.04062222318&0.08008008008&0.11817051419&0.03995747037\\
			\hline
			200&0.02879561418&0.056056056056&0.093947284323&0.03374368125\\
			\hline
			300&0.023531040266&0.046046046046&0.082115930412&0.03052496056\\
			\hline
			500&0.0182391885&0.036036036036&0.06928996488&0.026816971166\\
			\hline
			1000&0.012903494&0.024024024024&0.05501378899&0.022341822538\\
			\hline
		\end{tabular}
	\end{center}
	\item 
	\[
	f(x)=x^3+x^2+1
	\]
	\begin{center}
		\begin{tabular}{ |p{1cm}|p{2.7cm}|p{2.8cm}|p{2.7cm}|p{2.8cm}|  }
			\hline
			&&&& \\
			n& $\mu_n(p=1)$&$\omega(f,\mu_n)(p=1)$&$\mu_{n1}$&$\omega_{2,1}(f,\mu_{n1})$\\
			\hline
			10& 0.123091490979&0.552776666221&0.2474488015&0.226192537158\\
			\hline
			100&0.04062222318&0.19385157335&0.11817051419&0.09274336714\\
			\hline
			200&0.02879561418&0.1370198767&0.093947284323&0.07145628374\\
			\hline
			300&0.023531040266&0.113007080332&0.082115930412&0.061485978306\\
			\hline
			500&0.0182391885&0.08879734372&0.06928996488&0.050993785272\\
			\hline
			1000&0.012903494&0.05948463952&0.05501378899&0.039701675432\\
			\hline
		\end{tabular}
	\end{center}
\end{enumerate}
\begin{remark}
	It is evident that the estimates in Theorem \ref{swetits} are better than Theorem \ref{quan1} for $2$ and $3$. However, we remark that the estimates obtained \ref{quan1} are applicable to more general Banach function spaces. For instance, in the Weighted $L_p$ space we get better estimates for our choice of weight function as given below.
\end{remark}
\begin{example}{\textbf{Weighted Lebesgue space}}
	Let $L_{p,w}(0,1)$, $1<p<+\infty$ be the weighted Lebesgue space of the collection of measurable functions $f$ on $[0,1]$ with norm
	\begin{equation*}
		\|f\|_{p,w}=(\int_{0}^{1}|f(x)|^pw(x)dx)^{\frac{1}{p}}
	\end{equation*}
	and the weight function $w$ satisfying the Muckenhoupt condition :
	\begin{equation*}
		\sup\limits_{E\subset [0,1]}\frac{1}{|E|}\int_{E}w(x)dx(\frac{1}{|E|}\int_{E}w(x)^{-\frac{1}{p-1}}dx)^{p-1}<+\infty
	\end{equation*}
	We take the sequence of positive Kantorovich polynomials $\{K_n\}$ on $X^S$ . We have $K_n(1)=1$. Now
	\begin{equation*}
		\|1\|_{p,w}=(\int_{0}^{1}w(x)dx)^{\frac{1}{p}}
	\end{equation*}
	We get $c=(\int_{0}^{1}w(x)dx)^{\frac{1}{p}}<+\infty$ by the Muckenhoupt condition.\\
	Let $f\in C[0,1]$. Hence Theorem \ref{quan1} gives 
	\begin{equation*}
		\|K_n(f)-f\|_{p,w}\leq((\int_{0}^{1}w(x)dx)^{\frac{1}{p}}+1)\omega(f,\mu_n),\ n=1,2,...
	\end{equation*}
	where
	\begin{equation*}
		\mu_n^2=\|K_n((x-.)^2)(x)\|_{p,w}
	\end{equation*}
\end{example}

Let $w(x)=\frac{1}{1+x^2}$, for $x\in [0,1]$ be the weight function. Below, we give a numerical interpretation of this quantitative form for $p=1$ and $p=2$ for some choice of functions in $C[0,1]$.
\begin{enumerate} 
	\item 
	\[
	f(x)=x^2
	\]
	\begin{center}
		\begin{tabular}{ |p{1cm}|p{2.5cm}|p{2.8cm}|p{2.5cm}|p{2.8cm}|  }
			\hline
			&&&&\\
			$n$&$\mu_n(p=1)$&$\omega(f,\mu_n)(p=1)$& $\mu_n(p=2)$ &$\omega(f,\mu_n)(p=2)$\\
			\hline
			10&0.109461& 0.2068706199&0.119883&0.225290147\\
			\hline
			100&0.036145&0.0707178843&0.040073&0.078415363\\
			\hline
			200& 0.025623&0.0505546198&0.028429&0.0560044795\\
			\hline
			300& 0.020939& 0.0411755085& 0.023238&0.045870826\\
			\hline
			500& 0.016230&0.032143936 &0.018016&0.0356830718\\
			\hline
			1000& 0.011482&0.0226745489 &0.012748&0.02504621748\\
			\hline
		\end{tabular}
	\end{center}
	\item 
	\[
	f(x) = \begin{cases} 
		2x & \text{if } 0 \leq x \leq 0.5 \\
		-2x + 2 & \text{if } 0.5 < x \leq 1 
	\end{cases}
	\]
	\begin{center}
		\begin{tabular}{ |p{1cm}|p{2.5cm}|p{2.8cm}|p{2.5cm}|p{2.8cm}|  }
			\hline
			&&&&\\
			$n$&$\mu_n(p=1)$&$\omega(f,\mu_n)(p=1)$& $\mu_n(p=2)$ &$\omega(f,\mu_n)(p=2)$\\
			\hline
			10&0.109461& 0.2188437687&0.119883&0.239647929\\
			\hline
			100&0.036145&0.0720144028&0.040073&0.0800160032\\
			\hline
			200& 0.025623&0.051210242&0.028429&0.056811362\\
			\hline
			300& 0.020939& 0.0416083216& 0.023238&0.046409281\\
			\hline
			500& 0.016230&0.0324064812 &0.018016&0.0360072014\\
			\hline
			1000& 0.011482&0.0228045609& 0.012748&0.025205041\\
			\hline
		\end{tabular}
	\end{center}
	\item 
	\[
	f(x)=x^3+x^2+1
	\]
	\begin{center}
		\begin{tabular}{ |p{1cm}|p{2.5cm}|p{2.8cm}|p{2.5cm}|p{2.8cm}|  }
			\hline
			&&&&\\
			$n$&$\mu_n(p=1)$&$\omega(f,\mu_n)(p=1)$& $\mu_n(p=2)$ &$\omega(f,\mu_n)(p=2)$\\
			\hline
			10&0.109461& 0.500526951&0.119883&0.5434091\\
			\hline
			100&0.036145&0.1748966169&0.040073&0.193701485\\
			\hline
			200& 0.025623&0.1254199035&0.028429&0.13882379\\
			\hline
			300& 0.020939& 0.102298556&0.023238&0.113881877\\
			\hline
			500& 0.01623&0.079970277 &0.018016&0.08872732\\
			\hline
			1000& 0.011482&0.0564928367 &0.012748&0.06237931\\
			\hline
		\end{tabular}
	\end{center}
\end{enumerate}
Let $w(x)=\frac{e^{-x^2}}{1+64x^2}$, for $x\in [0,1]$.
\begin{enumerate} 
	\item 
	\[
	f(x)=x^2
	\]
	\begin{center}
		\begin{tabular}{ |p{1cm}|p{2.5cm}|p{2.8cm}|p{2.5cm}|p{2.8cm}|  }
			\hline
			&&&&\\
			$n$&$\mu_n(p=1)$&$\omega(f,\mu_n)(p=1)$& $\mu_n(p=2)$ &$\omega(f,\mu_n)(p=2)$\\
			\hline
			10&0.047490& 0.0925713048&0.076281&0.1466217228\\
			\hline
			100&0.015330&0.030174948&0.025363&0.0497747879\\
			\hline
			200&0.010852&0.021487634&0.017990&0.035290154\\
			\hline
			300& 0.008864&0.0175260497&0.014704&0.028992595\\
			\hline
			500& 0.006868&0.013556462 &0.011400&0.02227899\\
			\hline
			1000&0.004857&0.0095788711 &0.008066&0.015939175\\
			\hline
		\end{tabular}
	\end{center}
	\item 
	\[
	f(x) = \begin{cases} 
		2x & \text{if } 0 \leq x \leq 0.5 \\
		-2x + 2 & \text{if } 0.5 < x \leq 1 
	\end{cases}
	\]
	\begin{center}
		\begin{tabular}{ |p{1cm}|p{2.5cm}|p{2.8cm}|p{2.5cm}|p{2.8cm}|  }
			\hline
			&&&&\\
			$n$&$\mu_n(p=1)$&$\omega(f,\mu_n)(p=1)$& $\mu_n(p=2)$ &$\omega(f,\mu_n)(p=2)$\\
			\hline
			10&0.047490& 0.094818963&0.076281&0.152430486\\
			\hline
			100&0.015330&0.030406081&0.025363&0.050410082\\
			\hline
			200&0.010852&0.02160432&0.017990&0.0356071214\\
			\hline
			300& 0.008864& 0.01760352&0.014704&0.029205841\\
			\hline
			500& 0.006868&0.01360272&0.011400&0.02240448\\
			\hline
			1000&0.004857&0.00960192 &0.008066&0.0160032\\
			\hline
		\end{tabular}
	\end{center}
	\item 
	\[
	f(x)=x^3+x^2+1
	\]
	\begin{center}
		\begin{tabular}{ |p{1cm}|p{2.5cm}|p{2.8cm}|p{2.5cm}|p{2.8cm}|  }
			\hline
			&&&&\\
			$n$&$\mu_n(p=1)$&$\omega(f,\mu_n)(p=1)$& $\mu_n(p=2)$ &$\omega(f,\mu_n)(p=2)$\\
			\hline
			10&0.047490& 0.22816333&0.076281&0.358283878\\
			\hline
			100&0.015330&0.075094187&0.025363&0.123500041\\
			\hline
			200&0.010852&0.0535453159&0.017990&0.087755579\\
			\hline
			300& 0.008864& 0.043699599&0.014704&0.0721647357\\
			\hline
			500& 0.006868&0.0338220819 &0.011400&0.0555106472\\
			\hline
			1000&0.004857&0.0239127147 &0.008066&0.0397524114\\
			\hline
		\end{tabular}
	\end{center}
\end{enumerate}
\begin{example}{\textbf{Grand Lebesgue space}}
	Let $X=L_{* p}(0,1)$, $1<p<+\infty$ be a grand-Lebesgue space of measurable functions $f$ on $[0,1]$ with the norm 
	\begin{equation*}
		\|f\|_{* p}=\sup\limits_{0<\epsilon<p-1}\epsilon^{\frac{1}{p-\epsilon}}\|f\|_{p-\epsilon}
	\end{equation*}
	Then $X^S=\{f:\lim\limits_{\epsilon\rightarrow+0}\epsilon\int_{0}^{1}|f(x)|^{p-\epsilon}dx=0\}$. Consider the Kantorovich polynomials $\{K_n\}$ on $X^S$. Here $c=p-1$ so that application of Theorem \ref{quan1} gives for $f\in C[0,1]$,
	\begin{equation*}
		\|K_n(f)-f\|_{* p}\leq p\omega(f,\mu_n),\ n=1,2,...
	\end{equation*}
	where
	\begin{equation*}
		\mu_n^2=\|K_n((x-.)^2)(x)\|_{* p}
	\end{equation*}
\end{example}
In this example 
\begin{equation*}
	\|f\|_{* p}\leq (p-1)\|f\|_1
\end{equation*}
Thus for $f\in C[0,1]$,
\begin{equation*}
	\|K_n(f)-f\|_{* p}\leq p\omega(f,\sqrt{\frac{(p-1)}{(n+1)}}),\ n=1,2,...
\end{equation*}

\section*{Concluding Remarks}
In Theorems \ref{quan1} and \ref{quan2}, we derived quantitative estimates for functions in $C[0,1]$, where the concept of the modulus of continuity applies. We have observed that better estimates can be obtained in the case of weighted Lebesgue spaces. Extending the study further to obtain sharper estimates is a promising direction. Introducing new moduli of smoothness for general Banach function spaces and obtaining sharp estimates, as demonstrated in \cite{swetits}, is an interesting problem we plan to consider in the future.
\section*{Acknowlegements}
V.B. Kiran Kumar is supported by the KSYSA-Research Grant by KSCSTE, Kerala. P.C. Vinaya is thankful to the University Grants Commission(UGC), India for financial support.

\nocite{*}
\bibliographystyle{amsplain}
\bibliography{References}

\end{document}